\documentclass{amsart}
\usepackage{amssymb}

\usepackage{graphicx}
\usepackage[dvips]{color}
\usepackage{xcolor}

\newcommand{\la}{\langle}
\newcommand{\ra}{\rangle}

\newcommand{\cF}{\mathcal{F}}

\newcommand{\llangle}{\langle \kern -0.2em \langle}	
\newcommand{\rrangle}{\rangle \kern -0.2em \rangle}
\newcommand{\inner}[1]{\left\langle  #1 \right\rangle }

\newcommand{\iinner}[1]{\llangle  #1 \rrangle }

\newtheorem{theorem}{Theorem}[section]

\newtheorem{proposition}[theorem]{Proposition}
\newtheorem{corollary}[theorem]{Corollary}

\theoremstyle{definition}
\newtheorem{definition}[theorem]{Definition}
\newtheorem{example}[theorem]{Example}

\theoremstyle{remark}
\newtheorem{remark}[theorem]{Remark}

\numberwithin{equation}{section}

\begin{document}

\title[Fractional Supershifts]{Fractional Supershifts and their associated Cauchy Evolution problems}



\author[N. Alpay]{Natanael Alpay}
\address{(NA)
	Department of Mathematics\\ 
	University of California, Irvine,
	Irvine, CA 92697 \\
	USA}
\email{nalpay@uci.edu}



\date{}

\dedicatory{}

\commby{}

\begin{abstract}
In this work, we extend the notion of supershifts and superoscillation sequence to fractional Fock spaces based on Gelfond–Leontiev fractional derivatives. We first introduce the fractional supershifts seqeunce, and then discuss the associated evolution Cauchy problem with the fractional supershifts as initial condition.
\end{abstract}

\maketitle



\section{Introduction}

The phenomenon of superoscillations arises from quantum mechanics, and deals with functions that are band-limited signals, i.e. on suitable intervals, oscillate faster than their highest Fourier component.
The underlying physical phenomenon has been known for decades in areas such as weak values in quantum mechanics \cite{Aharonov1988}. 
With time, the study of superoscillations extended to \emph{supershifts}, introduced to study how superoscillatory data evolve under fundamental equations such as Schr\"odinger, Dirac, and Klein--Gordon. From a mathematical viewpoint, superoscillations arise when bounded spectra of ``small'' Fourier components combine to produce effective shifts far exceeding the spectral support, and have since led to new generating functions and connections to classical special functions, including Bernstein and Hermite families. \cite{ColomboKraussharSabadiniSimsek2023}.
For more mathematical theory and its applications to the Schr\"odinger evolution, see \cite{AharonovEtAl2017Memoir, yakir}; extensions of the  definitions and results to several variables appear in \cite{AharonovEtAl2016SeveralVars,ColomboPintonSabadiniStruppa2023}. \\

Examples of Supershift associated to reproducing kernel Hilbert space (RKHS) kernel were first studied in \cite{AlpayDeMartinoDiki2025ComplexRepresenter}.
In the present work, we  extend this idea to Supershift associated to the weight function of fractional Fock spaces. \\



Fractional derivatives extend the classical notion of derivative to noninteger orders,
the theory is motivated from modern signal and image processing, where one looks for frames obtained by discretizing integral transforms (e.g., wavelet and Gabor systems). Rooted in quantum mechanics and information theory (von Neumann, Gabor), the objective is to represent signals by time–frequency atoms with near-minimal joint spread. Practically, continuous transforms must be replaced by discrete, stable, and redundant systems; this raises lattice-density conditions in parameter space to guarantee a frame.
In the classical Gabor case, the work of Gr\"ochenig and Lyubarskii connects Gabor frames with Hermite windows to the standard orthonormal basis of the Fock space via the Bargmann transform, yielding concrete lattice criteria \cite{GL2007}.
Within this viewpoint, operator framework using Gelfond–Leontiev (GL) type derivatives, places fractional operators in a Hilbert-space setting and relate them to sampling and interpolation questions \cite{AlpayCerejeirasKahler2022}; where a more general class of Gabor-type frames is used. While the standard Fock framework is tied to classical differentiation and multiplication operators,  \cite{AlpayCerejeirasKahler2022} replace them with Gelfond–Leontiev derivatives (see \cite{K1994} for further details). The Gelfond–Leontiev derivatives encompasses many important operators as special cases, including Caputo and Riemann–Liouville fractional derivatives and difference–differential operators associated with finite reflection groups (Dunkl operators). The former have wide applications—from fractional mechanics to grey-noise models in stochastic processes—while the latter arise naturally in the analysis of Calogero–Sutherland–Moser $n$-particle systems. 
In this paper we adopt the GL perspective, and work in the \textit{fractional Fock space} based on GL derivatives.\\









Section 1 provides an introduction to fractional derivatives and superoscillation. Section 2 presents the necessary preliminaries, in particular, the construction of the fractional Fock space related to the Gelfond–Leontiev operators.
In Section 3, we define fractional superoscillations.
Sections 4 and 5 address the Cauchy problem associated with the fractional superoscillation sequence taken as the initial condition: Section 4 analyzes oscillatory integrals and distributional limits, which are subsequently used in Section 5 to find the solution.

\section{Preliminaries}

To study supershifts in the fractional Fock spaces related to Gelfond-Leontiev operators, we follow the buildup of generalized differentiation from \cite{AlpayCerejeirasKahler2022}.

\subsection{Generalized fractional derivatives}


Fractional derivatives via the Gelfond-Leontiev operators, depends on the type of function one takes. Choosing functions on the disk would yield the Hardy space, choosing entire functions leads to the Fock space, which is the space of interest in the paper. 


\begin{definition}
	Let
	\begin{equation} \label{Eq:EntireFunction}
		\varphi(z) =\sum_{k=0}^\infty \varphi_k z^k,
	\end{equation}
	be an entire function of order $\rho >0$ and degree $\sigma > 0,$ that is  $\lim_{k \rightarrow \infty} k^{\frac{1}{\rho}} \sqrt[k]{|\varphi_k|} =\left(\sigma e \rho\right)^{\frac{1}{\rho}}.$ We define the Gelfond-Leontiev (GL) operator of generalized differentiation with respect to $\varphi,$ denoted as $D_\varphi,$
	as the operator acting on an analytic function $f(z) =\sum_{k=0}^\infty a_k z^k,  |z|<1,$ as
	\begin{equation} \label{Eq:p1+p2}
		f(z) =\sum_{k=0}^\infty a_k z^k \quad \mapsto \quad D_\varphi f(z) =\sum_{k=1}^\infty a_k \frac{\varphi_{k-1}}{\varphi_k} ~z^{k-1}. \end{equation}
\end{definition}

Hence, under the condition on $\varphi$ that $\lim \sup_{k \rightarrow \infty} \sqrt[k]{\left|\frac{\varphi_{k-1}}{\varphi_k}\right|} =1$ by the Cauchy-Hadamard formula, we see that the series in (\ref{Eq:p1+p2})
inherit the same radius of convergence $R>0$ of the original series $f.$
Also, we would like to point out that the function $\varphi$ acts as a replacement of the \textit{exponential function} for the Gelfond-Leontiev operator of generalized differentiation and thus satisfies the identity $D_\varphi \varphi =\varphi.$\\

For a given entire function $\varphi$, we associate an analytic weight function $K_\varphi$ (on an appropriate domain) using the following relation.
\begin{equation}\label{eq:varn}
	\frac{1}{\varphi_n} = \mathcal{M}(K_\varphi(- \cdot))(n+1),
\end{equation}
where $ \mathcal{M}$ denotes the Mellin transform.
This relation induces a discrete reproducing kernel given by
\begin{equation}
	\label{eq: k_varphi}
	k_\varphi(n, k) := \varphi_n \delta_{n,k}, \qquad n, k \in \mathbb N_0,
\end{equation}	
with the corresponding discrete reproducing kernel Hilbert space is defined as 
\begin{equation}
	\label{eq: H(k_varphi)}
	\mathcal H(k_\varphi) := \Big\{ \underline f := (f_n)_{n=0}^\infty :  \| \underline f \|_{\ell^2_\varphi}^2  = \sum_{n=0}^\infty \frac{|f_n|^2}{\varphi_n} < \infty  \Big\},
\end{equation}
(sometime denoted by $\ell^2_\varphi$), and so for all $f\in \mathcal H(k_\varphi)$ we have the reproducing property
\begin{equation}
	\label{eq:f_in_H}
	f(z)=\la k_\varphi(z, \cdot), f \ra_{\ell^2_\varphi}.
\end{equation}
For two entire functions $f(z) =\sum_{k=0}^\infty f_k z^k, g(z) = \sum_{k=0}^\infty g_k z^k,$
the space $\mathcal{H}(k_\varphi)$ is equipped with the inner product
\[
\langle f,g\rangle_{2,{\varphi}}
:= \sum_{n=0}^\infty \frac{f_n\,\overline{g_n}}{\varphi_n}.
\]
The fractional Fock space $\mathcal{F}_\varphi$, is defined as 
\begin{equation}
	\label{eq: H(K_varphi)}
	\cF_{\varphi} := \Big\{ f := \sum_{n=0}^\infty f_n z^n :  \iinner{f, f}_{\cF_\varphi}  < \infty  \Big\},	
\end{equation}
with respect to the weighted measure $d\mu(w) = K_\varphi(-|w|^2) dxdy$, and endowed with the weighted inner product \begin{equation}\label{Eq:InnerProduct_Fock}
	\iinner{f, g}_{\cF_\varphi} = \frac{1}{\pi} \int_{\mathbb{C}} \overline{f(z)} g(z) K_\varphi(-|z|^2) dxdy,
\end{equation}
The space $\cF_\varphi$ is a fractional RKHS with the continuous kernel 
$$k_\varphi(z, w) := \varphi(\bar{z} w) = \sum_{n=0}^\infty \varphi_n\cdot (\bar{z} w)^n
\quad z,w\in\mathbb{C}
.$$

In fact, there is an isometry between $\ell^2_\varphi$ and $\cF_{\varphi}$
between $f(z)= \sum_{k=0}^\infty f_k z^k$ in $\cF_\varphi$ and its sequence of coefficients $(f_k)_{k=0}^\infty$ in $\ell^2_\varphi$
. i.e., these weighted inner products  should be related by the identity $\inner{f, g}_{\ell^2_\varphi} = \iinner{f, g}_{\cF_\varphi}$, that is 
\begin{equation}\label{Eq:Relation}
	\inner{f, g}_{\ell^2_\varphi} = \sum_{k=0}^\infty  \frac{\overline{f_k} g_k}{ \varphi_k} = \frac{1}{\pi} \int_{\mathbb{C}} \overline{f(z)} g(z) K_\varphi(-|z|^2) dxdy = \iinner{f, g}_{\cF_\varphi}.
\end{equation}

\begin{remark}
	The function $K_\varphi(z)$ represents the weight function of the Fractional Fock space, and the function $k_\varphi(z,w)=\varphi(z\bar w)$ is it's associated reproducing kernel.
\end{remark}
\begin{remark}
	Equation \eqref{eq:varn} reduces the determination of the measure $d\mu(z)=K_\varphi (-|z|^2) dx dy$ to either a problem of inversion of the Mellin transform or to a Stieltjes moment problem, for which  a sufficient condition for the determination of $K_\varphi$ consists in the Carleman's condition
	\begin{equation}\label{eq:carleman}
		\sum_{n=1}^\infty\tfrac{1}{\varphi_n^{\frac{1}{2n}}}=+\infty.
	\end{equation}
\end{remark}

\begin{example}\label{exponential} 
	Given the weight function $K_\varphi(z) = e^z$, from the inverse Melin transform with $\tilde K_{\varphi}(x) =  K_\varphi(-x) $, we get the kernel coefficients 
	$$
	\frac{1}{\varphi_n} 
	=
	\mathcal{M}(e^{-x}) (n+1)
	=
	\int_0^\infty  x^n  e^{-x}  dx 
	= n!
	. $$
	To show that the weight satisfies the moment problem, we verify \eqref{eq:carleman}, that is
	\begin{align*}
		\lim_{n\to\infty} \varphi_{2n}^{-1/(2n)}
		&=\lim_{n\to\infty} (2n!)^{\frac{1}{2n}} = \infty.
	\end{align*}
	Hence, the condition is satisfied.
	On a side note, starting with $\varphi(x)=e^{x}$, we remark another way to identify the weight function $K_\varphi(x) = e^{x}$, as 
	$$n! = \Gamma(n+1) = \int_0^\infty  x^n  e^{-x} dx.$$
\end{example}


\begin{example}\label{ex:mlf}
	Let $\varphi$ be the Mittag-Leffler function defined as:
	\begin{equation}\label{Eq:ML_function}
		E_{\frac{1}{\rho},\mu}(z) =\sum_{k=0}^\infty \tfrac{z^k}{\Gamma\left(\mu +\frac{k}{\rho}\right)}, \qquad \rho>0, ~\mu \in \mathbb{C}, ~{\rm Re}(\mu)>0,
	\end{equation}
	with coefficients $\varphi_k  =\frac{1}{\Gamma\left(\mu +\frac{k}{\rho}\right)}.$ 
	To verify \eqref{eq:carleman}, we first remark the following asymptotic formula for the Gamma function {\cite[A.1.24]{ml-book}}, 
	\[
	\Gamma(z) \sim \sqrt{\frac{2\pi}{z}} \left(\frac{z}{e}\right)^{z} \left[ 1 + \frac{1}{12 z} + \frac{1}{288z^2} +\cdots\right],
	\]
	as $z\to \infty$ with $|\text{arg}\;z|<\pi$ (a consequence of Stirling’s formula).
	Since we consider real $z$ as it approches infinity, and $\Gamma(z+1)=z\Gamma(z)$, we get 
	\[
	\lim_{x\to\infty} \Gamma(x+1) =\lim_{x\to\infty} \sqrt{2\pi x} \left(\frac{x}{e}\right)^{x}.
	\]
	To show that the weight satisfies the moment problem, we verify \eqref{eq:carleman}, 
	\begin{align*}
		\lim_{k\to\infty} \varphi_{2k}^{-1/(2k)} 
		&= \lim_{k\to\infty} \left(  \sqrt{2\pi k} \left(\frac{k}{e}\right)^{k} \right)^{1/(2k)}=\infty.
	\end{align*}
	Hence, the condition \eqref{eq:carleman} is satisfied. On the other hand, recall that the $\Gamma$ function is defined by a Mellin transform integral
	\[
	\Gamma(qn+1 )  = \int_{0}^{\infty} t^{qn} e^{-t}   
	\]
	Letting $r^2 = t^q$, yields $r^{2/q} = t$ and $\frac{2}{q} r^{\frac{2}{q} - 1} dr = dt$, so that
	\[
	\int_0^\infty r^{2n+1} d\tilde{\mu}(r) = \int_0^\infty r^{2n+1} \cdot \left( \frac{1}{q\pi} r^{-\frac{2}{q}} e^{-r^{2/q}} \right) dr.
	\]
	Thus, the measure functions  given by
	\[
	K_\varphi(-|z|^2)  = \frac{1}{q\pi} |z|^{\frac{2}{q} - 2} e^{-|z|^{\frac{2}{q}}}.
	\]
	This gives the Mittag-Leffler reproducing kernel Hilbert space denoted by $M L^2(\mathbb{C} ; q)$, first introduced by Rosenfeld in \cite{ross}, and in \cite{MLB, thermal}, the associated Bargman transform and the Thermal states were studied, respectively.
\end{example}



In \cite{AlpayCerejeirasKahler2022}, the generalized Bargmann transform $\tilde{\mathcal{B}}: L^2(\mathbb{R}) \to \mathcal{F}_\varphi$ given by
\begin{align*}
	\tilde{\mathcal{B}}f(z)&=\int_\mathbb{R}  \left( \sum_{n=0}^{\infty} \overline{h_n(x)} \sqrt{\varphi_n}z^n \right)  f(x)  dx
	=\int_{\mathbb{R}} A(z,x) f(x),
\end{align*}
where $A(z,x)$ is the Bargmann Kernel.
We observe that $\tilde{\mathcal{B}}$ is unitay by construction.
This gives the action of the raising and lowering operators $a^*$ and $a$ on the Hermite functions as
\begin{align}\label{eq:14}
	a^*h_{n-1} & = \sqrt{\frac{\varphi_{n-1}}{\varphi_{n}}}h_n , \qquad
	a h_{n}  = \sqrt{\frac{\varphi_{n-1}}{\varphi_{n}}}h_{n-1}, \qquad n=1,2,\ldots 
\end{align}
Let us emphasize that $\mathcal{B}$ acts now as an intertwining operator in the following way
\begin{align}\label{eq:15}
	\mathcal{B} af &= 
	D_\varphi \mathcal{B} f , \qquad
	\mathcal{B} a^*f=z\mathcal{B} f.
\end{align}

\begin{theorem}
	The system $\big\{e_n\big\}_{n\ge 0}$ with $e_n:=\sqrt{\varphi_n}\,z^n$ is an ONB of $\cF_{\varphi}$.
\end{theorem}
\begin{proof}
	First we show orthonormality.
	For $m,n\ge 0$ we have
	\[
	\langle e_n,e_m\rangle_{\cF_{\varphi}}
	=\sqrt{\varphi_n} \sqrt{\varphi_m} \underbrace{\langle z_n,z_m\rangle_{\mathcal{F}_\varphi}}_{=\frac{\delta_{n,m}}{\varphi_n}}=1.
	\]
	Thus $\{e_n\}$ is orthonormal.
	Next, we show the reproducing property,
	for fixed a $w$, the kernel section $k_w(z):=k(z,w)=\sum_{n\ge 0}\varphi_n z^n \overline{w}^{\,n}$ belongs to $\cF_{\varphi}$ since
	\[
	\|k_w\|_{\cF_{\varphi}}^2
	=\sum_{n=0}^\infty \frac{|\varphi_n \overline{w}^{\,n}|^2}{\varphi_n}
	=\sum_{n=0}^\infty \varphi_n |w|^{2n}
	=\varphi(|w|^2)<\infty.
	\]
	(Since $\varphi$ is entire, we do not need to worry about the analyticity of $\varphi(|z|^2$), however, if $\varphi$ is holomorphic in a disk or radius $r$, we would need to restrict $|z|\in {r,\sqrt{r}})$.
	See
	\[
	\langle f, k_w\rangle_{\cF_{\varphi}}
	=\sum_{n=0}^\infty \frac{f_n\,\overline{\varphi_n \overline{w}^{\,n}}}{\varphi_n}
	=\sum_{n=0}^\infty f_n w^n
	=f(w),\quad\text{for}\quad f(z)=\sum_{n\ge 0} f_n z^n\in\cF_{\varphi}
	\]
	so $k_\omega$ is indeed the reproducing kernel of $\cF_{\varphi}$.
	Lastly, we show  Completeness. 
	\[
	\langle f ,e_n\rangle_{\cF_{\varphi}}
	= \sum_{m=0}^{\infty} f_m\langle z^m, e_n\rangle_{\cF_{\varphi}}
	= \sum_{m=0}^{\infty} f_m \sqrt{\varphi_n}\langle z^m, z^n\rangle_{\cF_{\varphi}}
	=\sum_{m=0}^{\infty} f_m \sqrt{\varphi_n}
	\frac{\delta_{n,m}}{\varphi_n}
	=\frac{f_n}{\sqrt{\varphi_n}}.
	\]
	So for all $f\in\cF_{\varphi}$ we have 
	\[
	f(z)=\sum_{n=0}^\infty f_n z^n
	=\sum_{n=0}^\infty \Big(\tfrac{f_n}{\sqrt{\varphi_n}}\Big)\, e_n
	=\sum_{n=0}^\infty 
	\langle f ,e_n\rangle_{\cF_{\varphi}}
	\,e_n,
	\qquad\text{and}\qquad
	\sum_{n=0}^\infty \Big|\tfrac{f_n}{\sqrt{\varphi_n}}\Big|^2 < \infty.
	\]

	So the closed span of $\{e_n\}$ is all of $\cF_{\varphi}$. Together with orthonormality, $\{e_n\}$ is an orthonormal basis.
	Finally, the we notice that indeed kernel decomposes as
	\[
	k(z,w)=\sum_{n=0}^\infty e_n(z)\,\overline{e_n(w)}
	=\sum_{n=0}^\infty \varphi_n z^n \overline{w}^{\,n}=\varphi(\overline{w}z),
	\]
\end{proof}


\begin{example}
	If we consider $\varphi(z) = E_q(z)$ corresponding to the space $M L^2(\mathbb{C} ; q)$ as in Example \ref{ex:mlf}, we retrieve back the MLB transform $B_q$ from \cite{ross} and the ONB,
	$$
	\left\{g_n(z)\right\}_{n=0}^{\infty}=\left\{\tfrac{z^n}{\sqrt{\Gamma(q n+1)}}\right\}_{n=0}^{\infty}.
	$$
\end{example}

\subsection{Superoscillations}
Superoscillations are sequences of functions that exhibit high-frequency behavior, despite being composed of only low-frequency components. Originally introduced in the context of weak measurements in quantum mechanics by Aharonov et al.~\cite{AharonovEtAl2017Memoir}, these sequences possess surprising analytic and asymptotic properties that have found recent application in functional analysis and mathematical physics, and has significant implications for numerical methods in signal processing, quantum mechanics, and wave propagation.\\



%

The phenomenon of superoscillations has been studied in \cite{AharonovEtAl2017Memoir} using the sequence of functions defined by
\begin{equation*}
	F_n(a,t) = \left[ \cos \frac{t}{n} + i a \sin \frac{t}{n} \right]^n,\quad n\in\mathbb{N},\;a>1.
\end{equation*}

For large $n$ (equivalently small $t$), we have the supershift property  
$$\lim_{n\to\infty} F_n(a,t)= e^{iat}.$$
when $ a = 1 $, this expression is exact for any $ n $. For $ a > 1 $, the function $ F_n(a,t) $ oscillates at a frequency higher than its highest Fourier component; this will be more evident after the following proposition.

\begin{proposition}[{\cite[Theorem 3.1.17]{AharonovEtAl2017Memoir}}]\label{prop:superosc}
	The superoscillation sequence can be written using the Fock space kernel $B(z,w) = e^{z\bar{w}}$ as follows
	\begin{equation}\label{eq:suposc}
		F_n(z,a)=\sum_{j=0}^n C_j(n, a) e^{iz \left(1 - \frac{2j}{n}\right)}=\sum_{j=0}^{n} C_j(n,a)B(z,z_j),
	\end{equation}
	where for $ j=0,\dots n$ we have
	\begin{equation}\label{eq:cjzj}
		C_j(n,a) = {n\choose j} \left(\frac{1+a}{2}\right)^{n-j} \left(\frac{1-a}{2}\right)^j \quad \text{and}\quad 
		z_j  = -i \left(1-\frac{2j}{n}\right).
	\end{equation}
\end{proposition}


\begin{corollary}
	The superoscillation sequence can be written using the Fock space weight function $K_\varphi(x)=e^{x}$, as 
	\begin{equation*}
		F_n(z,a)=\sum_{j=0}^n C_j(n, a) e^{iz \left(1 - \frac{2j}{n}\right)}=\sum_{j=0}^{n} C_j(n,a) K_\varphi(-zz_j),
	\end{equation*}
	where the constants $C_j(n,a)$ and $z_j$ are given by \eqref{eq:cjzj} 
\end{corollary}

\begin{remark}
	Equation \eqref{eq:suposc} is where the superosilatory behavior is the most apparent, as the frequency $1-\frac{2j}{n}$ of each Fourier component of $F_n(z,a)$ is bounded by $1$, but at the limit $e^{iat}$, the frequency $a$ can be made arbitrary big.
\end{remark}

\begin{proposition}
	\label{prop:lim}
	Let $F_n(z,a)$ be the function defined in \eqref{eq:suposc}. Then for any $z$ in some  compact sets of $\mathbb{C}$, we have 
	\[
	\lim_{n\to\infty} F^{(k)}_n(z,a) = (ia)^ke^{iza}.
	\]
	In particular 
	\[\lim_{n\to\infty} F_n^{(k)}(0,a) = (ia)^k.\]
\end{proposition}

\begin{proof}
	
	Following proof of {\cite[Proposition 3.1.16.]{AharonovEtAl2017Memoir}}, 
	let $F_n(z,a)$ as defined in \eqref{eq:suposc}, we directly compute the derivatives of the function $F_n(z, a)=g_n^n(z)$ where
	$$
	g_n(z)=\cos \left(\frac{z}{n}\right)+i a \sin \left(\frac{z}{n}\right).
	$$
	We see that
	$$
	F_n^{\prime}(z, a)=n g_n^{n-1}(z) g_n^{\prime}(z)=n \frac{g_n^{\prime}(z)}{g_n(z)} F_n(z, a).
	$$
	
	We are now interested in the asymptotic behavior of $F_n^{\prime}$ 
	when $n$ goes to infinity and $z$ is in a compact subset of $\mathbb{C}$. 
	From {\cite[theorem 3.1.8]{AharonovEtAl2017Memoir}}, we know $F_n(z, a)\to e^{iaz}$ uniformly as $n\to\infty$ on every compact subset of $\mathbb{R}$, but can be extended naturally to $\mathbb{C}$.
	So it is enough to compute the asymptotic behavior of $g_n^{\prime} / g_n$ we have ,
	$$
	\frac{g_n^{\prime}(z)}{g_n(z)}=\frac{1}{n} \frac{-\sin \left(\frac{z}{n}\right)+i a \cos \left(\frac{z}{n}\right)}{\cos \left(\frac{z}{n}\right)+i a \sin \left(\frac{z}{n}\right)} \sim \frac{1}{n} i a.
	$$
	So
	\begin{align*}
		\lim_{n\to\infty} F'_n(z,a) &= \lim_{n\to\infty} n \frac{g'_n}{g_n} F_n(z,a)
		=\lim_{n\to\infty} ia F_n(x,a)
		=iae^{iax}.
	\end{align*}
	Repeating the process with the $k^{th}$ derivative, we get 
	\[
	\lim_{n\to\infty} F^{(k)}(z,a) = (ia)^k e^{iza}.
	\]
	Taking $z=0$ implies the result.
\end{proof}

\begin{corollary}\label{corollary:ia}
	Let $F_n(z,a)$ be the function defined in \eqref{eq:suposc}. Then
	\[
	\lim_{n\to\infty } \sum_{j=0}^{n} C_j(n,a) i^k\left( 1-\frac{2j}{n} \right)^k 
	=(ia)^k.
	\]
	In particular, for $k=0$ 
	\[
	\lim_{n\to\infty} \sum_{j=0}^{n} C_j(n,a) = 1.
	\]
	By {\cite[Remark 3.1.11]{AharonovEtAl2017Memoir}} this equality holds even more generally as
	$$F(0,a) = \sum_{j=0}^{n} C_j(n,a)=  1, $$  
\end{corollary}

We now have all the necessary background and tools to extend the notion of the superoscillations sequence to the generalized fractional Fock space setting.

	\section{Fractional Supershifts}

	Proposition \ref{prop:superosc} was the motivation of the present work in which the superoscillation sequence was written in terms of the weight function of the Fock space, i.e $K_\varphi(x)=e^{x}$. In this work we would want to extend the idea to generalized supershifts written in terms of the fractional Fock space kernel.\\
	
	
	Let us remind the fractional Fock space $\mathcal{F}_\varphi$ introduced in \cite{AlpayCerejeirasKahler2022} based on Gelfond-leontiev derivatives, defined as 
	\begin{equation}
		\label{eq: H(K_varphi)}
		\cF_{\varphi} := \Big\{ f := \sum_{n=0}^\infty f_n z^n :  \iinner{f, f}_{\cF_\varphi}  < \infty  \Big\},	
	\end{equation}
	with respect to the weighted measure $d\mu(w) = K_\varphi(-|w|^2) dxdy$, and endowed with the weighted inner product \begin{equation}\label{Eq:InnerProduct_Fock}
		\iinner{f, g}_{\cF_\varphi} = \frac{1}{\pi} \int_{\mathbb{C}} \overline{f(z)} g(z) K_\varphi(-|z|^2) dxdy,
	\end{equation}
	with the associated reproducing kernel
	\begin{equation}\label{eq:rk}
		k_\varphi(z,w) = \varphi(z\bar{w}) =\sum_{n=0}^{\infty} \varphi_n (z\bar{w} )^n.    
	\end{equation}

	Let $K_\varphi$ be the associated measure to $\varphi$ via \eqref{eq:varn}.
	We now define the analogs of the superoscillations, called generalized supershifts, or fractional supershifts, in a similar spirit to Proposition \ref{prop:superosc} in which we saw that classical superoscillation written in terms of $K_\varphi(x)=e^{x}$ corresponding to $\varphi(x)=e^{x}$ .

	\begin{definition}\label{def:FracSuper}
		For a fractional Fock space $\mathcal{F}_\varphi$ with weight function $K_\varphi$ anaytic on $|z|<r$, for $n\geq 1$ define the fractional supershift sequence as
		\begin{equation}\label{eq:frac-osc}
			F_{n,\varphi}(z,a)=\sum_{j=0}^{n} C_j(n,a) K_\varphi(-z_j\cdot z). 
		\end{equation}
		for all $|z|<r$, $a>1$, and $C_j(n,a),\;z_j$ as given in \eqref{eq:cjzj}.
	\end{definition}

	To show that this definition indeed defines a sequence holding the supershift property as in Proposition \ref{prop:lim}, we give the following proposition.

	\begin{proposition}[Supershift property]\label{prop:supshift}
		For $a>1$, $z_j$ and $C_j(n,a)$ defined by \eqref{eq:cjzj}.
		Using a weight function $K_\varphi$ analytic on $|z|<r$, the function $F_{\varphi,n}$ extends the notion of superoscillations via the supershift property.
		\[
		\lim_{n\to\infty} F_{\varphi,n}(z,a) =  K_\varphi(-iza),\quad \forall\;|z|<r.
		\]
	\end{proposition}
	\begin{proof}

		For a fixed $n>0$, and $z_j$, $C_j(n,a)$ given by \eqref{eq:cjzj}.
		Since $K_\varphi$ analytic on $|z|<r$ consider the power series expansion 
		\[K_\varphi(z) = \sum_{k=0}^{\infty} a_k z^k, \quad |z|<r. \]
		The superoscillation sequence $F_{n,\varphi}$ can then be expended as 
		\begin{align*}
			F_{n,\varphi}(z,q)
			=\sum_{j=0}^{n}C_j(n,a)K_\varphi\left(iz\left(1-\frac{2j}{n}\right)\right)
			=
			\sum_{k=0}^{\infty} a_k\cdot  (iz)^k  \left(
			\sum_{j=0}^{n}C_j(n,a) \left(1-\frac{2j}{n}\right)^k\right).
		\end{align*}
		We want to show 
		\[
		\lim_{n\to\infty} |F_{n,\varphi}(z,a)- K_\varphi(-iza)|=0.
		\]
		Since $K_\varphi$ is analytic on $|z|<r$, we have
		\[
		|F_{n,\varphi}(z,a)- K_\varphi(-iza)|\leq \sum_{k=0}^{\infty}  |a_k|  |z|^k \left|
		\sum_{j=0}^{n} C_j(n,a) i^k\left( 1-\frac{2j}{n} \right)^k
		-(ia)^k
		\right|.
		\]
		By Proposition \ref{prop:lim} it follows
		\[
		\lim_{n\to\infty } \sum_{j=0}^{n} C_j(n,a) i^k\left( 1-\frac{2j}{n} \right)^k
		=\lim_{n\to\infty} F_{n}^{(k)}(0,a)=(ia)^k.
		\]
		Hence $\lim_{n\to\infty} |F_{n,\varphi}(z,a)-K_\varphi(-iza)|=0.$ and the result follows.
	\end{proof}

	Therefore, $F_{n,\varphi}$ indeed defines a superoscillation; the next question would be to study the moments and Cauchy evolution problems in the next sections.

	\begin{remark}
		For $\varphi (z)= e^{z}$, we have $K_\varphi(z)=e^z$, and so we retrieve the classical superoscillation sequence.
	\end{remark}

	\begin{example}\label{ex:ml}
		For $q\geq 0$ Considering $\varphi(z) =  E_q(z)$, from Example \eqref{ex:mlf}, we know 
		$$K_\varphi(-|z|^2)  = \frac{1}{q\pi} |z|^{\frac{2}{q} - 2} e^{-|z|^{\frac{2}{q}}}.$$
		So, for $a\geq 1$ we get the fractional supershifts associated with the Mittag-Leffler Fock space defined as 
		\[
		F_{n,q}(z,a)=\frac{1}{q\pi} \sum_{j=0}^{n} C_j(n,a)
		|zz_j|^{\frac{1}{q} - 1} e^{-|zz_j|^{\frac{1}{q}}}
		\xrightarrow{n \to \infty } 
		\frac{1}{q\pi} |az|^{\frac{1}{q} - 1} e^{-a|z|^{\frac{1}{q}}}
		\]
	\end{example}

	\begin{remark}
		Another natrual way to define generalized supershifts is using the kernel. Eamples of such supershifts were first studied in \cite{AlpayDeMartinoDiki2025ComplexRepresenter} as a solution to the least-square minimization problem.. Let $k_\varphi = \varphi$ be the kernel of RKHS $\mathcal{F}_\varphi$, define
		\[
		\tilde F_{n,\varphi}(z,a)=\sum_{j=0}^{n} C_j(n,a)\varphi(z, z_j)\quad \forall z\in\mathbb{C}. 
		\]
		These sequences would also satisfy the supershift property
		\[
		\lim_{n\to\infty} \tilde F_{\varphi,n}(z,a) =  \varphi(-iza),\quad \forall z\in\mathbb{C}.
		\]
		using a similar proof as in Proposition \ref{prop:supshift} and the expansion $ \varphi(z) = \sum_{k=0}^{\infty}\varphi_k z^k.$
		Moreover, this definition depends only on $\varphi$ and the property holds for all $z\in\mathbb{C}$.
		However, for such a generalized supershifts sequence, we did not find a direct way to solve the moments and Cauchy evolution problem (see \eqref{eq:where_we_need_relation}), as they do not possess the Melin transform relation \eqref{eq:varn}.
	\end{remark}

	\begin{remark}
		In the work we consider an entire function $\varphi$, however the theory can be extended to analytic $\varphi$, as in {\cite[Example 3.4]{AlpayCerejeirasKahler2022}} in which we have 
		\[
		\varphi(z) = \sum_{k=0}^\infty \frac{z^k}{\pi\cot(\pi(k+1))}
		\xrightarrow{\text{associated weight}} 
		K_\varphi (z)  = \frac{1}{1-z}.
		\]
		analytic on the unit disk.
		Moreover, even when $\varphi$ is entire, the associated weight function $K_\varphi$ obtained via the inverse Mellin transform need not be entire. Using 
		\[
		\mathcal{M}(t^{\alpha}e^{-t})(s) = \int_{0}^{\infty} t^{s-1} t^{\alpha} e^{-t} dt = \Gamma(s+\alpha).
		\]
		for $\alpha=\frac{1}{2}$, we get the weight function 
		\[
		\mathcal{M}[K_\varphi(-\cdot))](s) = i\Gamma\left(s+\frac12\right) = \frac{1}{\varphi_{s-1}}.
		\]
		Thus, one can construct an explicit example where $\varphi$ is entire but $K_\varphi$ is not,
		\[
		\varphi(z) = \sum_{k=0}^\infty \frac{z^n}{i\Gamma\left(k+\frac32\right) }
		\xrightarrow{\text{associated weight}} \;
		K_\varphi (z)  =  \sqrt{z} e^{-z} .
		\]
	\end{remark}

	\section{Oscillatory Integrals and Distributional Limits}

	To study the evolution problem, we must first study the notions of limits of integrals in distributions. One can skip to remark \ref{rmk:alternate} at the end of this section for an alternative way; however, the content of this section is interesting. 
	Let $\mathcal{S}(\mathbb{R})$ denote the Schwartz space of rapidly decreasing smooth functions.
	The theory of Fourier integral operators was first developed by H\"ormander in \cite{hormanderFIO1}, and further extended by Duistermaat \cite{duistermaat}, Treves \cite{treves}, and Taylor \cite{taylorPDE3}.
	A Fourier integral operator $ T $ acting on a function $ u \in \mathcal{S}(\mathbb{R}) $ is defined by the formula
	\[
	Tu(x) = \int_{\mathbb{R}} \int_{\mathbb{R}} e^{i \phi(x,\xi) - i y \xi} a(x,\xi) u(y) \, dy \, d\xi,
	\]
	where $ \phi(x,\xi) $ is a smooth real-valued phase function satisfying certain non-degeneracy conditions, and $ a(x,\xi) $ is a symbol of order $ m \in \mathbb{R} $.
	In one dimension, the Schwartz kernel of an FIO 
	with phase $\phi(x,\xi)$
	can be written as the oscillatory integral
	\[
	K(x,y)=\int_{\mathbb{R}} e^{\,i(\phi(x,\xi)-y\xi)}\,a(x,\xi)\,d\xi,
	\]
	and the operator acts by
	\[
	(Tu)(x)=\int_{\mathbb{R}} K(x,y)\,u(y)\,dy.
	\]
	
	\begin{remark}
		Applying $T$ to the \emph{constant function} $u\equiv 1$ gives
		\[
		(T1)(x)=\int_{\mathbb{R}}\!\!\left(\int_{\mathbb{R}}
		e^{-iy\xi}\,dy\right)a(x,\xi)e^{i\phi(x,\xi)}\,d\xi
		=2\pi\,a(x,0)\,e^{i\phi(x,0)}
		\]
		(up to Fourier normalization). Thus $(T1)(x)$ is {not}
		$\int e^{i\phi(x,\xi)}a(x,\xi)\,d\xi$. The latter appears when the kernel
		is evaluated at $y=0$, i.e. for $T\delta_0$
		\[
		(T\delta_0)(x)=K(x,0)=\int_{\mathbb{R}} e^{i\phi(x,\xi)}\,a(x,\xi)\,d\xi.
		\]
	\end{remark}

	In order to solve the Schr\"odinger evolution problem associated with the generalized superoscillation in the next section, we must first study the regularized oscillatory integral 
	\[
	I_m^\epsilon(x,t) = \int_{\mathbb{R}} \lambda^m e^{-(i t + \epsilon)\lambda^2 + i \lambda x} \, d\lambda
	\]
	as the kernel of a Fourier integral operator (FIO). We then show that its distributional limit as $ \epsilon \to 0^+ $ fits into the standard framework of FIOs.
	The regularized kernel $ I_m^\epsilon(x,t) $ can be written in the form
	\[
	I_m^\epsilon(x,t) = \int_{\mathbb{R}} a(\lambda) e^{i \phi(x,\lambda)} e^{-\epsilon \lambda^2} \, d\lambda,
	\]
	with real-valued, smooth, non-degenerate  phase function $ \phi(x,\lambda) = \lambda x - t \lambda^2 $
	satisfying
	$\frac{\partial^2 \phi}{\partial x \, \partial \lambda} = 1 \ne 0,
	$
	and  amplitude $ a(\lambda) = \lambda^m $.
	This integral is well-defined for all $x \in \mathbb{R}$ and $t \in \mathbb{R}$ due to the negative real part of the exponent ensuring integrability, i.e. since 
	\[
	\left| \lambda^m e^{-(i t + \epsilon)\lambda^2 + i \lambda x} \right| = |\lambda|^m e^{-\epsilon \lambda^2},
	\]
	and the function $\lambda \mapsto |\lambda|^m e^{-\epsilon \lambda^2}$ belonging to $L^1(\mathbb{R})$ for every fixed $\epsilon > 0$, we have absolute convergence of the integral.\\

	Let 
	$\mathcal{S}'(\mathbb{R})$ denote the dual of $\mathcal{S}(\mathbb{R})$ ,  the space of tempered distributions.
	We have the following propositions.

	\begin{proposition}
		Let $ m \in \mathbb{N}_0 $, $ t \in \mathbb{R} $, and $ \epsilon > 0 $. Then the function
		\[
		I_m^\epsilon(x,t) = \int_{\mathbb{R}} \lambda^m e^{-(i t + \epsilon)\lambda^2 + i \lambda x} \, d\lambda
		\]
		is the kernel of a Fourier integral operator associated with the phase function $ \phi(x,\lambda) = \lambda x - t \lambda^2 $ and amplitude $ a(\lambda) = \lambda^m $.
		Moreover, the limit
		\[
		I_m(x,t) := \lim_{\epsilon \to 0^+} I_m^\epsilon(x,t)
		\]
		exists in $ \mathcal{S}'(\mathbb{R}_x) $ and defines a tempered distribution that is also the kernel of a Fourier integral operator in the sense of distribution theory.
	\end{proposition}

	%

	\begin{proof}
		The regularized integral $ I_m^\epsilon(x,t) $ fits the classical form of a Fourier integral operator 
		\begin{align*}
			I_m^\epsilon(x,t)
			&= \int_{\mathbb{R}} \lambda^m e^{-(it+\epsilon)\lambda^2 + i\lambda x}\,d\lambda= \int_{\mathbb{R}} a(\lambda)\,e^{i\phi(x,\lambda)}\,e^{-\epsilon\lambda^2}\,d\lambda.
		\end{align*}
		with phase $ \phi(x,\lambda) = \lambda x - t \lambda^2 $, which satisfies the standard non-degeneracy conditions
		\[
		\frac{\partial^2 \phi}{\partial x \, \partial \lambda} = \frac{\partial}{\partial x} \left( \frac{\partial \phi}{\partial \lambda} \right) = \frac{\partial}{\partial x}(x - 2t\lambda) = 1.
		\]

		The amplitude $ \lambda^m \in S^m $ is a classical symbol of order $ m $.
		The function $a(\lambda) = \lambda^m$ is a smooth function on $\mathbb{R}$, and it satisfies for any multi-index $\alpha$,
		\[
		|\partial_\lambda^\alpha a(\lambda)| \leq C_\alpha (1 + |\lambda|)^{m - |\alpha|},
		\]
		which means $a \in S^m(\mathbb{R})$, the class of classical symbols of order $m$, so $a(\lambda)$ is an admissible amplitude for an FIO.
		Now we will show convergence for a fixed $\epsilon > 0$, the integrand
		\[
		f_\epsilon(x,\lambda) := \lambda^m e^{-(i t + \epsilon)\lambda^2 + i \lambda x},
		\]
		is smooth and rapidly decaying in $\lambda$, 
		as we have the pointwise bound
		\[
		|f_\epsilon(x, \lambda)| = |\lambda|^m e^{-\epsilon \lambda^2},
		\]
		which is integrable over $\mathbb{R}$ for each $m \in \mathbb{N}_0$. Therefore, $I_m^\epsilon(x,t)$ is a smooth function of $x$ for each fixed $\epsilon > 0$.\\
		
		Let $\varphi \in \mathcal{S}(\mathbb{R}_x)$ be a Schwartz test function. We consider the pairing of $I_m^\epsilon(x,t)$ with $\varphi(x)$
		\[
		\langle I_m^\epsilon(\cdot,t), \varphi \rangle = \int_{\mathbb{R}} \left( \int_{\mathbb{R}} \lambda^m e^{-(i t + \epsilon)\lambda^2 + i \lambda x} \, d\lambda \right) \varphi(x) \, dx.
		\]
		Since we have absolute integrability of $I_m^{\epsilon}$, 
		we apply Fubini's theorem to exchange the order of integration
		\[
		\langle I_m^\epsilon(\cdot,t), \varphi \rangle = \int_{\mathbb{R}} \lambda^m e^{-(i t + \epsilon)\lambda^2} \left( \int_{\mathbb{R}} e^{i \lambda x} \varphi(x) \, dx \right) d\lambda = \int_{\mathbb{R}} \lambda^m e^{-(i t + \epsilon)\lambda^2} \widehat{\varphi}(-\lambda) \, d\lambda,
		\]
		where $\widehat{\varphi}$ denotes the Fourier transform of $\varphi$, which also lies in $\mathcal{S}(\mathbb{R})$.
		We now consider the limit $\epsilon \to 0^+$. 
		The pointwise limit is,
		\[
		\lim_{\epsilon\to 0^+} \lambda^m e^{-(i t + \epsilon)\lambda^2} \widehat{\varphi}(-\lambda) =
		\lambda^m e^{-i t \lambda^2} \widehat{\varphi}(-\lambda),
		\]
		For all $\lambda \in \mathbb{R}$, since $\widehat{\varphi}$ is a Schwartz function, it satisfies rapid decay, for any $N>0$, we have
		\[
		|\widehat{\varphi}(-\lambda)| \leq C_N(1+|\lambda|)^{-N}.
		\]
		Taking $N = m+2$, we obtain
		\[
		|\lambda^m e^{-(i t+\epsilon)\lambda^2}\widehat{\varphi}(-\lambda)| \leq |\lambda|^m|\widehat{\varphi}(-\lambda)| \leq C_{m+2}(1+|\lambda|)^{-2},
		\]
		which provides an integrable dominating function. 
		Thus, by the Dominated Convergence Theorem, we can interchange the limit and integral to conclude
		\[
		\lim_{\epsilon \to 0^+}\langle I_m^\epsilon(\cdot,t),\varphi\rangle = \int_{\mathbb{R}}\lambda^m e^{-i t\lambda^2}\widehat{\varphi}(-\lambda)\,d\lambda.
		\]
		This defines a continuous linear functional on $\mathcal{S}(\mathbb{R})$, and therefore, the limit
		\[
		I_m(x,t) := \lim_{\epsilon \to 0^+} I_m^\epsilon(x,t)
		\]
		exists in the sense of tempered distributions.
		This completes the proof.
	\end{proof}

	This limit would allow us to give a more precise formulation of the integral $I_m(x,t)$ as the limit of $I^{\epsilon}_{m}(x,t)$ for which we can find the exact value using the following formula [from $9.253$ (cite integral book)] as follows,
	\[
	\int_{-\infty}^{\infty} (ix)^v e^{-\beta^2 x^2 -iqx}\;dx=
	2^{-\frac{v}{2}} \sqrt{\pi} \beta^{-v-1} \exp\left( \frac{-q^2}{8\beta^2} \right) D_v\left(\frac{q}{\beta \sqrt{2}}\right),
	\]
	for $\text{Re}\beta>0,\;\text{Re}(v)>-1,\; \text{arg}\;ix = \frac{\pi}{2}\text{sign}\;x$;
	and 
	\[
	D_n(z) = 2^{-n/2} e^{-z^2/4} H_n\left( \frac{z}{\sqrt{2}} \right)
	\]
	where $H_n$ are the $n$'th Hermite polynomials.
	Since $x\in \mathbb{R}$, it holds $\text{arg}\;ix = \frac{\pi}{2}\text{sign}\;x$. For $\epsilon>0$, taking $v = m$, $\beta^2 = (it+\epsilon) $, $ q = -x$,
	we get
	\begin{align*}
		I_m^\epsilon(x,t)&=\int_{\mathbb{R}}\lambda^m e^{ -(it + \epsilon) \lambda^2 + i \lambda x } \; d\lambda \\
		&= \frac{1}{i^m} 2^{-\frac{m}{2}} \sqrt{\pi} (it+\epsilon)^{-\frac{m+1}{2}} \exp\left( \frac{-x^2}{ 8(it+\epsilon) } \right) D_m \left( \frac{-x}{\sqrt{2(it+\epsilon)}} \right)\\
		&= \frac{1}{i^m} 2^{-{m}} \sqrt{\pi} (it+\epsilon)^{-\frac{m+1}{2}} \exp\left( \frac{-x^2}{ 4(it+\epsilon) } \right)H_m\left(\frac{-x}{2\sqrt{it+\epsilon}} \right).
	\end{align*}
	Then we see that the integral $I_m(x,t)$ that we need to evaluate is equivalent to 
	\begin{align*}
		I_m(x,t) = \lim_{\epsilon\to 0^+}I_m^\epsilon(x,t)
		&= \lim_{\epsilon\to 0^+} 	\int_{\mathbb{R}}\lambda^m e^{ -(it + \epsilon) \lambda^2 + i \lambda x } \; d\lambda \\
		&= \lim_{\epsilon\to 0^+}\frac{1}{i^m} 2^{-{m}} \sqrt{\pi} (it+\epsilon)^{-\frac{m+1}{2}} \exp\left( \frac{-x^2}{ 4(it+\epsilon) } \right)
		H_m\left(\frac{-x}{2\sqrt{it+\epsilon}} \right).
	\end{align*}
	To show that the limit is well defined, we need to show that the limit exists. 
	Since the complex-valued Hermite polynomials are entire and coincide with the real Hermite polynomials \cite{Gors,DunklXu2014,Folland1989},
	we get
	\[
	\lim_{\epsilon\to 0^+} H_m\left(\frac{-x}{2\sqrt{it+\epsilon}} \right)=H_m\left(\frac{-x}{2\sqrt{2t}}(1-i) \right).
	\]
	And so we can find the following closed form for $I_m$,
	\begin{equation}
		\label{eq:Im}
		I_m(x,t) = 
		\frac{1}{i^m} 2^{-{m}} \sqrt{\pi} (it)^{-\frac{m+1}{2}}
		\exp\left( \frac{-x^2}{ i4t } \right)
		H_m\left(\frac{-x}{2\sqrt{it}} \right).
	\end{equation}
	

	\begin{remark}\label{rmk:alternate}
		From {\cite[Theorem 6.5]{AlpayDeMartinoDiki2025ComplexRepresenter}}, we have a general formula which stats that for $a,b\in\mathbb{C}\setminus\{0\}$ with $\Re(a)>0$. Then for $n\in\mathbb{N}$ we have
		\begin{equation}
			\label{eq:thm6.5}
			\int_{-\infty}^{\infty} x^{n} e^{-a x^{2}+b x}\,dx
			=
			\sqrt{\frac{\pi}{a}}\,
			e^{\frac{b^{2}}{4a}}\,
			(-i)^{n}
			\left(\frac{1}{2\sqrt{a}}\right)^{n}
			H_{n}\!\left(i\,\frac{b}{2\sqrt{a}}\right).
		\end{equation}
		using 
		$
		a = it,\, b = ix,\, n=m
		$.
		Then
		\[
		\sqrt{\frac{\pi}{a}}=\sqrt{\pi}\,(it)^{-1/2},\qquad
		\exp\!\left(\frac{b^2}{4a}\right)
		=\exp\!\left(\frac{(ix)^2}{4it}\right)
		=\exp\!\left(\frac{-x^2}{i4t}\right),
		\]
		and
		\[
		(-i)^m\left(\frac{1}{2\sqrt{a}}\right)^m
		=(-i)^m\,2^{-m}(it)^{-m/2}
		=i^{-m}\,2^{-m}(it)^{-m/2}
		=\frac{1}{i^m}\,2^{-m}(it)^{-m/2}.
		\]
		Moreover, the Hermite argument becomes
		$
		i\frac{b}{2\sqrt{a}}
		=i\frac{ix}{2\sqrt{it}}
		=-\frac{x}{2\sqrt{it}}.
		$
		Substituting these into \eqref{eq:thm6.5} yields
		\[
		\int_{\mathbb{R}}\lambda^m e^{-it\lambda^2+i\lambda x}\,d\lambda
		=
		\frac{1}{i^m}2^{-m}\sqrt{\pi}\,(it)^{-\frac{m+1}{2}}
		\exp\!\left(\frac{-x^2}{i4t}\right)
		H_m\!\left(-\frac{x}{2\sqrt{it}}\right),
		\]
		which is exactly the integral we calculated.
		However, we calculated this theorem after this section was written; moreover, I find this section's content to be interesting enough to keep.
	\end{remark}

		\section{ Evolution Problems: 
			Schr\"odinger equation for the free particle}
		
		A solution to Schr\"odinger evolution of superoscillations and supershifts was given by  Aharonov et al in \cite{yakir}.
		For 
		$$H\psi(x,t):=-\frac{\partial^2\psi(x,t)}{\partial x^2}.$$
		They solved the evolution in time of the following Cauchy problem
		\begin{equation*}
			\begin{cases}
				i  \frac{\partial \psi(x,t)}{\partial t} =H \psi(x,t),\\
				\psi(x,0) = F_n(x,a).
			\end{cases},
		\end{equation*} 
		In \cite{AlpayDeMartinoDiki2025ComplexRepresenter} a solution to the $ L^2$-regularized Cauchy problem with initial condition
		\[        \psi(x,0) = \phi_n(x) F_n(x,a).
		\]
		for different $\phi_n(x)$ such as the Gaussian and Hermite polynomials.
		In this section, we want to solve the evolution problem
		
		\begin{equation}\label{eq:evo_problem}
			\begin{cases}
				i  \frac{\partial \psi(x,t)}{\partial t} =H \psi(x,t),\\
				\psi(x,0) = F_{n,\varphi}(x,a),
			\end{cases}
		\end{equation}
		where $F_{n,\varphi}$ is the supershift sequence associated to a valid weight function $K_\varphi$.
		\begin{theorem}
			The solution to Problem \ref{eq:evo_problem}, is given by 
			\[
			\psi(x,t) =
			\frac{1}{\sqrt{2\pi}} 
			\sum_{m=0}^{\infty}
			b_n(m)
			\exp\left( \frac{-x^2}{ i4t } \right)
			H_m\left(\frac{-x}{2\sqrt{it}} \right),
			\]
			where $b_n(m)$ is given by
			\[
			b_m = \frac{\sqrt{\pi} }{m!}
			\frac{1}{\varphi_{m+1}}
			2^{-{2}}  (it)^{-\frac{m+1}{2}}
			\sum_{j=1}^{n} C_j(n,a)z_j^{1-m},
			\]
			and $C_j,z_j$ defined by \eqref{eq:cjzj}.
		\end{theorem}
		
		\begin{proof}

			To solve the problem, we apply the Fourier transform on $i\frac{\partial}{\partial t} \psi(x,t) = - \frac{\partial^2}{\partial x^2} \psi(x,t)$ and then solve the resulting ODE. We get 
			\[ 
			i\frac{\partial}{\partial t } \hat{\psi}(\lambda,t) = \lambda^2 \hat{\psi}(\lambda,t)
			\implies 
			\hat{\psi}(\lambda,t) =  c(\lambda) e^{-i\lambda^2 t}  ,
			\]
			where the constant $c(\lambda)$ is calculated from the initial condition $\hat{\psi}(\lambda,0) = c(\lambda)$ for all $\lambda \in \mathbb{R}$ with $\psi(x,0)=F_{n,\varphi}(x,a)$.
			Using the initial condition and $ z_j=-i\left( 1-\frac{2j}{n} \right)$, we get
			\begin{equation}\label{eq:where_we_need_relation}
				\begin{aligned}
					\widehat{\psi}(\lambda,0) = c(\lambda) & =  \frac{1}{\sqrt{2\pi}}
					\int_{\mathbb{R}} e^{-i\lambda x}\psi(x,0)dx\\
					&= \frac{1}{\sqrt{2\pi}}\int_{\mathbb{R}} e^{-i\lambda x } F_{n,\varphi}(x,a)dx\\
					&= \frac{1}{\sqrt{2\pi}}\sum_{j=0}^nC_j(n,a) \int_{\mathbb{R}} e^{-ix \lambda } K_\varphi\left( x z_j \right)  dx\\
					&= \frac{1}{\sqrt{2\pi}}\sum_{m=0}^{\infty} \frac{(i\lambda)^{m}}{m!} \sum_{j=1}^{n} C_j(n,a) \int_{\mathbb{R}} x^m K_\varphi(x z_j)dx.
				\end{aligned}
			\end{equation}
			Making the change of variable $y = xz_j$, $dy = z_j dx$, and using relation \eqref{eq:varn}, we get 
			\begin{align*}
				\int_{\mathbb{R}} x^m K_\varphi(x z_j)dx & = z_j^{1-m} \int_{\mathbb{R}} y^m K_\varphi(y)dy
				= z_j^{1-m} \mathcal{M}(K_\varphi)(m+1)
				=  z_j^{1-m} \frac{1}{\varphi_{m+1}}.
			\end{align*}
			Combining the terms, it follows
			\[
			c(\lambda) = \sum_{m=0}^{\infty} \frac{(i\lambda)^{m}}{m!}
			\frac{1}{\varphi_{m+1}}
			\underbrace{\sum_{j=1}^{n} C_j(n,a)z_j^{1-m}}_{:=c_n(m)} 
			=
			\sum_{m=0}^{\infty} \frac{(i\lambda)^{m} c_n(m) }{m!}
			\frac{1}{\varphi_{m+1}}.
			\]
			Substituting $c(\lambda)$ back into $\hat\psi(\lambda,t)$, we get
			\[
			\hat{\psi}(\lambda,t) 
			=  \frac{1}{\sqrt{2\pi}}
			\sum_{m=0}^{\infty} \frac{(i)^{m} }{m!}
			\frac{1}{\varphi_{m+1}}c_n(m) \;\lambda^m e^{ -i\lambda^2t }.
			\]
			To find $\psi$, we are left to take the inverse Fourier transform, and we get
			\begin{align*}
				\psi(x,t) & = \frac{1}{\sqrt{2\pi}} \int_{\mathbb{R}} e^{i\lambda x} \hat{\psi}(\lambda , t ) d\lambda 
				=
				\frac{1}{\sqrt{2\pi}} 
				\sum_{m=0}^{\infty} \frac{(i)^{m} }{m!}
				\frac{1}{\varphi_{m+1}}c_n(m)
				\underbrace{\int_{\mathbb{R}}\lambda^m e^{ -i t \lambda^2 + i \lambda x } \; d\lambda}_{= I_m(x,t) }.
			\end{align*}
			In the previous, section we studied integral $I_m(x,y)$, using result \eqref{eq:Im}, we get
			\[
			I_m(x,t) = 
			\frac{1}{i^m} 2^{-{m}} \sqrt{\pi} (it)^{-m-1}
			\exp\left( \frac{-x^2}{ i4t } \right)
			H_m\left(\frac{-x}{2\sqrt{it}} \right).
			\]
			Substituting as follows completes the problem, 
			\begin{align*}
				\psi(x,t) & = \frac{1}{\sqrt{2\pi}} \int_{\mathbb{R}} e^{i\lambda x} \hat{\psi}(\lambda , t ) d\lambda \\
				&=
				\frac{1}{\sqrt{2\pi}} 
				\sum_{m=0}^{\infty} \frac{(i)^{m} }{m!}
				\frac{1}{\varphi_{m+1}}c_n(m)
				\left[
				\frac{1}{i^m} 2^{-m} \sqrt{\pi} (it)^{-\frac{m+1}{2}} \exp\left( \frac{-x^2}{ i4t } \right)
				H_m\left(\frac{-x}{2\sqrt{it}} \right)
				\right]\\
				&=
				\frac{1}{\sqrt{2\pi}} 
				\sum_{m=0}^{\infty} 
				\underbrace{\frac{\sqrt{\pi} }{m!}
					\frac{1}{\varphi_{m+1}}
					c_n(m)
					2^{-{m}}  (it)^{-\frac{m+1}{2}}}_{=b_n(m)}
				\exp\left( \frac{-x^2}{ i4t } \right)
				H_m\left(\frac{-x}{2\sqrt{it}} \right).
			\end{align*}

		\end{proof}

		\subsection*{Acknowledgments} The author would like to thank Professors Paula Cerejeiras and Uwe Kaehler for their valuable discussions and insightful suggestions. Their guidance, particularly in approaching the integral problem from a distributional point of view, was especially helpful.
		
		
		\bibliographystyle{amsplain}
		\bibliography{bibtex}
		
	\end{document}